\documentclass[11pt,a4paper]{article}

\usepackage[mathcal]{eucal}
\usepackage{color}
\usepackage{latexsym}
\usepackage{amsfonts}
\usepackage{amsmath}
\usepackage{amsthm}
\usepackage{amssymb}
\usepackage{mathrsfs}
\usepackage{graphicx}
\usepackage{algorithm}
\usepackage{algorithmic}

\newtheorem{theorem}{Theorem}[section]

\newlength{\programindent}
\newcommand{\BIGOP}[1]{\mathop{\mathchoice%
{\raise-0.22em\hbox{\huge $#1$}}%
{\raise-0.05em\hbox{\Large $#1$}}{\hbox{\large $#1$}}{#1}}}

\newcounter{noteCounter}

\begin{document}

% ----------------------------------- TITLE MATTER ----------------------------------------------

\title{Not All Saturated 3-Forests Are Tight}

\author{
% Heidi Gebauer\thanks{}\\
% \vspace{0.1cm} \\
% \small Institute for Theoretical Computer Science\\
% \small Department of Computer Science\\
% \small ETH Z\"urich, 8092 Z\"urich, Switzerland\\
% \small \texttt{gebauerh@inf.ethz.ch}
% \and
% Anna Gundert\thanks{supported by the Swiss National Foundation (SNF Project 200021-125309)}\\
% \vspace{0.1cm} \\
% \small Institute for Theoretical Computer Science\\
% \small Department of Computer Science\\
% \small ETH Z\"urich, 8092 Z\"urich, Switzerland\\
% \small \texttt{anna.gundert@inf.ethz.ch}
% \and
% Robin A. Moser\thanks{}\\ 
% \vspace{0.1cm} \\ 
% \small Institute for Theoretical Computer Science\\
% \small Department of Computer Science\\
% \small ETH Z\"urich, 8092 Z\"urich, Switzerland\\
% \small \texttt{robin.moser@inf.ethz.ch}
% \and
Heidi Gebauer,\; Anna Gundert\thanks{Supported by the Swiss National Foundation (SNF Project 200021-125309).},\; Robin A. Moser
 \vspace{0.1cm} \\
 \small Institute for Theoretical Computer Science\\
 \small Department of Computer Science\\
 \small ETH Z\"urich, 8092 Z\"urich, Switzerland\\
 \small \texttt{\{gebauerh, anna.gundert, robin.moser\}@inf.ethz.ch}
\vspace*{0.1cm} \\
\and
\setcounter{footnote}{5}
Yoshio Okamoto\thanks{Supported by Grant-in-Aid for Scientific Research from Ministry
  of Education, Science and Culture, Japan, and
  Japan Society for the Promotion of Science.}
\vspace*{0.1cm} \\
\small Center for Graduate Education Initiative\\
\small Japan Advanced Institute of Science and Technology\\
\small 1-1 Asahidai, Nomi, Ishikawa, 923-1292 Japan\\
\small \texttt{okamotoy@jaist.ac.jp}
}
\date{\today}

\maketitle

\begin{abstract}
A basic statement in graph theory is that every inclusion-maximal forest is connected, i.e.\ a tree.
Using a definiton for higher dimensional forests by Graham and Lov\'asz and the connectivity-related notion of tightness for hypergraphs introduced by Arocha, Bracho and Neumann-Lara in, we provide an example of a saturated, i.e.\ inclusion-maximal $3$-forest that is not tight. This resolves an open problem posed by Strausz.
%We provide an example of a saturated $3$-forest that is not tight, thus resolve an open problem by ....
\end{abstract}

\smallskip
\noindent
\section{Introduction}
In this note, we consider a generalization of the graph-theoretic concepts of trees and forests to $k$-uniform hypergraphs. Several approaches to this can be found in the literature. The notion for forests studied here was defined by Graham and Lov\'asz \cite{Lo}, the one for trees by Arocha, Bracho and Neumann-Lara \cite{ArBrNL}.
The basic approach taken here is to find higher-dimensional analogues of two graph properties characterizing trees: acyclicity and connectivity.

For graphs, being acyclic and connected is the same as being an inclusion-maximal acyclic graph, i.e.\ a tree. It was an open question whether this also holds for the generalization presented here. We present a counterexample, a $3$-uniform hypergraph that has the generalized acyclicity property and is inclusion-maximal but does not satisfy the higher-dimensional connectivity property.

For acyclic graphs, having exactly $n-1$ edges is the same as being an inclusion-maximal acyclic graph. For the generalizations presented here, this is not true: an inclusion-maximal hypergraph satisfying the generalized acyclicity property need not have the maximum number of edges. The  $3$-uniform hypergraph presented in this paper is an example for this.
Actually, for graphs any two out of the three properties acyclicity, connectivity and having exactly $n-1$ edges imply the third. Using the same generalizations for acyclicity and connectivity as in this note, Parekh \cite{Par} shows that this does not hold for these generalized properties either. He suggests a stricter higher-dimensional analogue of connectivity to work around this.

\section{Preliminaries }
In what follows, $k$-uniform hypergraphs will be called $k$-graphs. By a \emph{$t$-coloring} of a $k$-graph $H=(V,E)$ we mean a \emph{surjective} mapping $c\colon V\to [t]$.  A set of vertices in a $t$-colored $k$-graph is \emph{polychromatic} if the vertices are colored with different colors.  
A polychromatic edge is called a \emph{rainbow} edge.
%A \emph{rainbow} edge in a $t$-colored $k$-graph is an edge that is colored with $k$ pairwise different colors.

In the (2-)graph case, a forest is an acyclic $2$-graph. This is the same as saying that every edge $e$ is a cut edge, i.e., there is a partition of the vertex set into two non-empty sets such that $e$ is the only edge connecting the two sets. This is our point of view for generalizing the concept of forests:
 we call a $k$-graph $H=(V,E)$ a \emph{$k$-forest} if for every edge $e \in E$ there is a $k$-coloring of $H$ that has $e$ as its only rainbow edge, i.e.\ there is a $k$-coloring $c$ such that for every $e' \in E$ we have: $c(e') = [k] \Leftrightarrow e'=e$.
Lov\'asz proved the following.
\begin{theorem}[Lov\'asz \cite{Lo}]
 A $k$-forest with $n$ vertices has at most $\binom{n-1}{k-1}$ edges.
\end{theorem}
Parekh \cite{Par} gave a second proof of this theorem. 

For (2-)graphs, a tree is the same as a connected forest. We now present a concept corresponding to connectivity in higher dimensions that we will need to define a notion of trees for hypergraphs.
In \cite{ArBrNL} Arocha, Bracho and Neumann-Lara introduce the following invariant:
the \emph{heterochromatic number} $\mathfrak{hc}(H)$ of a $k$-graph $H=(V,E)$ is the minimum number $t$ of colors such that any $t$-coloring of $H$ has a rainbow edge:
%an edge that is colored with $k$ different colors:
\[\begin{split}
   \mathfrak{hc}(H)& = \min \{t \mid \forall \, t\text{-colorings}\: c\text{ of } H\;\,\exists \, e\!\in\! E\colon\, |c(e)| =k\}\\
& = \max \{t\;|\; \exists \, \text{a }\,t\text{-coloring}\: c\text{ of } H \text{ without a rainbow edge}\} + 1. 
  \end{split}
\]
Obviously, $k \leq \mathfrak{hc}(H) \leq n$. A $k$-graph $H$ is called \emph{tight} if $\mathfrak{hc}(H)=k$.

One can quickly see that a $2$-graph is tight if and only if it is connected: both notions describe the property that for any partition of the vertex set into two non-empty sets there is an edge connecting the two sets.
This inspires the definition of a \emph{$k$-tree} as a tight $k$-forest. Equivalently, one can define a $k$-tree to be a tight $k$-graph which does not remain tight after the removal of any edge. This is because an edge whose removal destroys tightness is the unique rainbow edge for some coloring and vice versa.

%The author of the post at the Open Problem Garden (Ricardo Strausz?) claims to be able to prove the following:

%\begin{theorem}
 %Every $k$-forest with $n$ vertices and $\binom{n-1}{k-1}$ edges is tight - and therefore a tree.
%\end{theorem}

A $k$-forest is called \emph{saturated} if no edge can be added to it without losing the property of being a $k$-forest --- one might also call it an ``inclusion-maximal'' $k$-forest.
Note that for $k>2$ there are saturated $k$-forests with $n$ vertices and less than $\binom{n-1}{k-1}$ edges. The $3$-graph presented in this paper is an example for this.

In the (2-)graph case we have that any saturated $2$-forest is tight: an inclusion-maximal forest is a tree. The question of whether this also holds for $k>2$ was posed by Strausz \cite{opg}.
We study the case $k=3$. While one can see that any saturated $3$-forest on $4$ or $5$ vertices is tight, this note presents a counterexample which shows:
\begin{theorem}
 There exists a $3$-forest on $6$ vertices that is saturated but not tight, i.e.\ not a $3$-tree.
\end{theorem}

\section{The Counterexample}
Our counterexample is the $3$-graph $H = (V,E)$ on the six vertices $V := \{1,2,3,4,5,6\}$, containing the following edges:
$$ E := \{ \{1,2,3\}, \{1,2,4\}, \{1,2,5\}, \{1,3,4\}, \{2,3,6\}, \{2,5,6\}, \{3,4,6\}, \{3,5,6\}  \}. $$ 
\begin{figure}
 \centering
 \includegraphics[scale=.7]{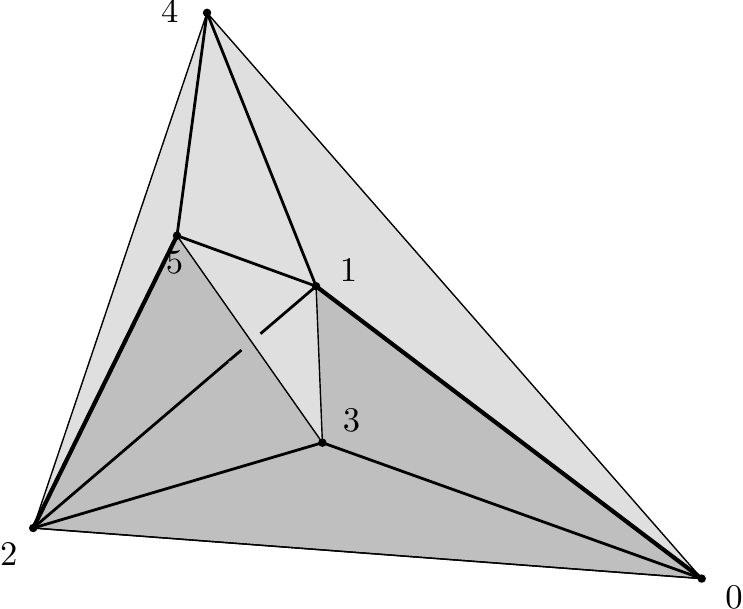}
 \caption{An illustration of the counterexample}\label{counterexample}
\end{figure}
See Figure \ref{counterexample} for an illustration.
%
%The example has been found using a simple computer program conducting all necessary brute-force explorations of the coloring space. Although such a program can be used to check the validity of the counterexample as well, let us give a short manually verifiable certificate.
To show that this is a counterexample we need to demonstrate that
\begin{itemize}
 \item[(i)] $H$ is a $3$-forest,
 \item[(ii)] $H$ is saturated and
 \item[(iii)] $H$ is not tight.
\end{itemize}

\subsection*{Proof of (i)}

To show (i), we need to exhibit a mapping $\gamma\colon E \to [3]^V$ such that for all $e, f \in E$, $\gamma(e)(f) = [3]$ if and only if $e=f$. (Remind that $[3]^V$ is the set of all maps $V\to [3]$.) The following table represents such a mapping.

\begin{center}
\begin{tabular}{|c|cccccc|}
\hline
$e$ & $\gamma(e)(1)$ & $\gamma(e)(2)$ & $\gamma(e)(3)$ & $\gamma(e)(4)$ & $\gamma(e)(5)$ & $\gamma(e)(6)$ \cr
\hline
$\{1,2,3\}$ & $3$ & $2$ & $1$ & $3$ & $2$ & $1$ \cr
$\{1,2,4\}$ & $1$ & $3$ & $1$ & $2$ & $1$ & $1$ \cr
$\{1,2,5\}$ & $3$ & $2$ & $2$ & $2$ & $1$ & $1$ \cr
$\{1,3,4\}$ & $3$ & $3$ & $1$ & $2$ & $1$ & $1$ \cr
$\{2,3,6\}$ & $3$ & $3$ & $2$ & $2$ & $1$ & $1$ \cr
$\{2,5,6\}$ & $3$ & $3$ & $1$ & $1$ & $2$ & $1$ \cr
$\{3,4,6\}$ & $3$ & $3$ & $3$ & $2$ & $1$ & $1$ \cr
$\{3,5,6\}$ & $1$ & $1$ & $3$ & $1$ & $2$ & $1$ \cr
\hline
\end{tabular}
\end{center}

For example,
\begin{align*}
\gamma(\{1,2,4\})(\{1,2,3\}) &= \{\gamma(\{1,2,4\})(1), \gamma(\{1,2,4\})(2), \gamma(\{1,2,4\})(3)\} = \{1,3\},\\
\gamma(\{1,2,4\})(\{1,2,4\}) &= \{\gamma(\{1,2,4\})(1), \gamma(\{1,2,4\})(2), \gamma(\{1,2,4\})(4)\} = \{1,2,3\}.
\end{align*}

It is easy to check that on each line, all polychromatic triples of vertices, except for the edge the coloring is being built for, are not contained in $E$.

\subsection*{Proof of (ii)}

To demonstrate (ii), we first reinspect the above table to see whether there are alternative ways to fill in certain rows. Formally, define $\Phi \colon E \to 2^{[3]^V}$
to map each $e \in E$ to the set of $3$-colorings in which $e$ is the only rainbow edge. For any $e \in E$, $\Phi(e)$ thus represents the set of colorings which could
be used to construct entry $\gamma(e)$ for proving (i). We want to explicitly determine certain values of $\Phi$. Moreover, let $\bar E := \binom{V}{3}\setminus E$ be the non-edges
in our example and let $\Delta \colon E \to 2^{\bar E}$ map each edge $e \in E$ to the set of all non-edges which are polychromatic under all colorings in $\Phi(e)$. We note
that any variant $H'$ that arises from $H$ by adding a new edge $e' \in \bar E$ cannot be a $3$-forest if there exists $e \in E$ such that $e' \in \Delta(e)$ because in such an $H'$, there is no $3$-coloring featuring $e$ as its only rainbow edge. We can thus rule out all these non-edges to demonstrate the saturation of $H$.

Let $c$ be a coloring that colors $\{1,2,3\}$ polychromatically, without loss of generality $c(1)=3, c(2)=2, c(3)=1$, and all other edges at most bichromatically. The presence of $\{1,2,4\} \in E$ 
requires $c(4) \in \{3,2\}$ and the presence of $\{1,3,4\}$ requires $c(4) \in \{3,1\}$, hence $c(4)=3$. The presence of the edge $\{2,3,6\}$ implies $c(6) \in \{2,1\}$ while
$\{3,4,6\}$ implies $c(6) \in \{1,3\}$ leaving no choice but $c(6)=1$. If $c(5)=3$, then $\{2,5,6\}$ features all colors. If $c(5)=1$, then $\{1,2,5\}$ does. That leaves no
other choice but $c(5)=2$. Therefore, the coloring used for $\gamma(\{1,2,3\})$ above is the unique member of $\Phi(\{1,2,3\})$ (up to relabelling of colors). This entails
that $\Delta(\{1,2,3\}) = \{\{1,2,6\},\{1,3,5\},\{1,5,6\},\{2,4,6\},\{3,4,5\},\{4,5,6\}\},$ since all of these edges are polychromatic under $\gamma(\{1,2,3\})$. None of these
non-edges can be added to $H$ without losing the $3$-forest property.

% THE FOLLOWING SECTION SEEMS TO BE SPARE
%Let now $c$ be a coloring that colors $\{1,2,4\}$ polychromatically. Without loss of generality, $c(1)=1, c(2)=3, c(4)=2$ as in $\gamma(\{1,2,4\})$. The presence of $\{1,2,3\} \in E$ requires
%$c(3) \in \{1,3\}$. The presence of $\{1,3,4\}$ on the other hand requires $c(3) \in \{1,2\}$ so that $c(3)=1$. The presence of both $\{2,3,6\}$ and $\{2,4,6\}$ enforces that
%$c(6) \in \{1,3\}$ and $c(6) \in \{2,3\}$, respectively, so that $c(6) = 3$. For $c(5)$ we have a choice between the two colors $1$ and $3$. $c(5) \ne 2$ is because otherwise,
%$\{2,5,6\}$ would become rainbow. We conclude that $\Phi(\{1,2,4\})$ consists of the two described colorings that differ only in $c(5)$ (and all colorings that arise from
%these by relabelling of the colors). Under these colorings, the non-edges $\Delta(\{1,2,4\}) = \{ \{2,3,4\}, \{2,4,6\} \}$ are always polychromatic.  None of these non-edges
%can thus be added to $H$ without losing the $3$-forest property.

Let now $c$ be a coloring that colors $\{1,3,4\}$ polychromatically. Without loss of generality, $c(1)=3, c(3)=1, c(4)=2$ as in $\gamma(\{1,3,4\})$. The presence of both $\{1,2,3\}$ and
$\{1,2,4\}$ requires via $c(2) \in \{1,3\}$ and $c(2) \in \{2,3\}$ that $c(2)=3$. Furthermore the presence of both $\{2,3,6\}$ and $\{3,4,6\}$ requires via $c(6) \in \{1,3\}$
and $c(6) \in \{1,2\}$ that $c(6)=1$. If $c(5)$ were $2$, then the edge $\{2,5,6\}$ would become rainbow. We are left with two options for $c(5)$, $1$ and $3$, which are
both possible, so that $\Phi(\{1,3,4\})$ consists of these two colorings and all their isomorphic variants. Under these colorings, the non-edges
$\Delta(\{1,3,4\}) = \{ \{1,4,6\}, \{2,3,4\}, \{2,4,6\} \}$ are always rainbow, such that none of these are eligible for addition to $H$ without producing a non-$3$-forest.

Let now $c$ be a coloring that colors $\{2,3,6\}$ polychromatically, without loss of generality $c(2)=3, c(3)=2, c(6)=1$ as in $\gamma(\{2,3,6\})$. Since there is $\{2,5,6\}$, we have
$c(5) \in \{1,3\}$  and since there is $\{3,5,6\}$, we also have $c(5) \in \{1,2\}$, leaving only $c(5)=1$. Since there are $\{1,2,3\}$ and also $\{1,2,5\}$, we conclude both $c(1) \in \{2,3\}$ and $c(1) \in \{1,3\}$, thus $c(1) = 3$. Finally, the presence of $\{1,3,4\}$ and $\{3,4,6\}$ requires $c(4) \in \{2,3\}$ and $c(4) \in \{1,2\}$
and thus $c(4)=2$. Up to relabelling, this is thus the only one option in $\Phi(\{2,3,6\})$.
%there is thus only this one option in $\Phi(\{2,3,6\})$.
This shows that $\Delta(\{2,3,6\}) = \{\{1,3,5\},\{1,3,6\},\{1,4,5\},\{1,4,6\},\{2,3,5\},\{2,4,5\},\{2,4,6\}\}$.

To be sure to have covered all necessary non-edges, the following table lists all triples in $\binom{V}{3}$.

%\begin{center}
%\begin{tabular}{|c|c|}
%\hline
%triple & contained in \cr
%\hline
%$\{1,2,3\}$ & $E$ \cr
%$\{1,2,4\}$ & $E$ \cr
%$\{1,2,5\}$ & $E$ \cr
%$\{1,2,6\}$ & $\Delta(\{1,2,3\})$ \cr
%$\{1,3,4\}$ & $E$ \cr
%$\{1,3,5\}$ & $\Delta(\{1,2,3\})$ \cr
%$\{1,3,6\}$ & $\Delta(\{2,3,6\})$ \cr
%$\{1,4,5\}$ & $\Delta(\{2,3,6\})$ \cr
%$\{1,4,6\}$ & $\Delta(\{1,3,4\})$ \cr
%$\{1,5,6\}$ & $\Delta(\{1,2,3\})$ \cr
%$\{2,3,4\}$ & $\Delta(\{1,3,4\})$ \cr
%$\{2,3,5\}$ & $\Delta(\{2,3,6\})$ \cr
%$\{2,3,6\}$ & $E$ \cr
%$\{2,4,5\}$ & $\Delta(\{2,3,6\})$ \cr
%$\{2,4,6\}$ & $\Delta(\{1,2,3\})$ \cr
%$\{2,5,6\}$ & $E$ \cr
%$\{3,4,5\}$ & $\Delta(\{1,2,3\})$ \cr
%$\{3,4,6\}$ & $E$ \cr
%$\{3,5,6\}$ & $E$ \cr
%$\{4,5,6\}$ & $\Delta(\{1,2,3\})$ \cr
%\hline
%\end{tabular}
%\end{center}

\begin{center}
\begin{tabular}{|c|c||c|c|}
\hline
triple & contained in & triple & contained in \cr
\hline
$\{1,2,3\}$ & $E$                 & $\{2,3,4\}$ & $\Delta(\{1,3,4\})$ \cr
$\{1,2,4\}$ & $E$                 & $\{2,3,5\}$ & $\Delta(\{2,3,6\})$ \cr
$\{1,2,5\}$ & $E$                 & $\{2,3,6\}$ & $E$ \cr
$\{1,2,6\}$ & $\Delta(\{1,2,3\})$ & $\{2,4,5\}$ & $\Delta(\{2,3,6\})$ \cr
$\{1,3,4\}$ & $E$                 & $\{2,4,6\}$ & $\Delta(\{1,2,3\})$ \cr
$\{1,3,5\}$ & $\Delta(\{1,2,3\})$ & $\{2,5,6\}$ & $E$ \cr
$\{1,3,6\}$ & $\Delta(\{2,3,6\})$ & $\{3,4,5\}$ & $\Delta(\{1,2,3\})$ \cr
$\{1,4,5\}$ & $\Delta(\{2,3,6\})$ & $\{3,4,6\}$ & $E$ \cr
$\{1,4,6\}$ & $\Delta(\{1,3,4\})$ & $\{3,5,6\}$ & $E$ \cr
$\{1,5,6\}$ & $\Delta(\{1,2,3\})$ & $\{4,5,6\}$ & $\Delta(\{1,2,3\})$ \cr
\hline
\end{tabular}
\end{center}

\subsection*{Proof of (iii)}

Last, in order to establish (iii), it suffices to give a single $3$-coloring that does not make any edge rainbow. In fact, if we use
the coloring $c$ defined as
\begin{center}
\begin{tabular}{|cccccc|}
\hline
$c(1)$ & $c(2)$ & $c(3)$ & $c(4)$ & $c(5)$ & $c(6)$ \cr
\hline
$3$ & $2$ & $2$ & $2$ & $2$ & $1$ \cr
\hline
\end{tabular}
\end{center}
then clearly no edge becomes rainbow.

\section*{Acknowledgements}
This work has been done at the 9th Gremo's Workshop on Open Problems (GWOP 2011).
The authors would like to thank Tobias Christ for producing Figure \ref{counterexample}, the organizers of GWOP and Emo Welzl for continuous support.

\bibliographystyle{amsalpha}
\bibliography{paper-arxiv-kforests11.bib}

\end{document}